\documentclass[11pt,leqno]{article} 

\usepackage{fullpage}

\usepackage{amssymb} 
\usepackage{amsmath} 

\date{}

\input xy 
\xyoption{all}
\title{Classification of simple linearly compact n-Lie superalgebras} 
\author{{\sc Nicoletta Cantarini}\thanks{Dipartimento di Matematica
Pura ed Applicata, Universit\`a di Padova, Padova, Italy
- Partially supported by Progetto di ateneo CPDA071244}
\and\setcounter{footnote}{6}
{\sc Victor G.\ Kac}\thanks{Department of Mathematics, MIT, Cambridge,
Massachusetts 02139, USA
- Partially supported by an NSF grant}}

\newtheorem{theorem}{Theorem}[section] 
\newtheorem{lemma}[theorem]{Lemma} 
\newtheorem{corollary}[theorem]{Corollary} 
\newtheorem{proposition}[theorem]{Proposition} 
\newtheorem{definition}[theorem]{Definition} 
\newtheorem{remark}[theorem]{Remark}

\def\Z{\mathbb{Z}} 
\def\ZZ{\mathbb{Z}}
\def\g{\mathfrak{g}}

\def\fg{\mathfrak{g}}

\def\F{\mathbb{F}}
\def\FF{\mathbb{F}}
\def\CC{\mathbb{C}}
\def\ZZ{\mathbb{Z}}

\def\0{\bar{0}}
\def\1{\bar{1}}
\def\O{{\mathcal O}}

\numberwithin{equation}{section} 

\newcommand{\st}[1]{\ensuremath{^{\scriptstyle \textrm{#1}}}}

\makeatletter

\def\enumerate{%
  \ifnum \@enumdepth >\thr@@\@toodeep\else
    \advance\@enumdepth\@ne
    \edef\@enumctr{enum\romannumeral\the\@enumdepth}%
      \list
        {\csname label\@enumctr\endcsname}%
        {\usecounter\@enumctr
          \addtolength{\leftmargin}{-\leftmargin}
          \settowidth{\labelwidth}{(99)}
          \itemindent = \labelwidth
          \addtolength{\itemindent}{\labelsep}
        \listparindent=1em      
          \def\makelabel##1{{##1}\hfill}
          }%
  \fi}

\makeatother

\newcommand{\alphaparenlist}{
  \renewcommand{\theenumi}{\alph{enumi}}%
  \renewcommand{\labelenumi}{(\theenumi)}%
}

\newcommand{\romanparenlistii}{
  \renewcommand{\theenumii}{\roman{enumii}}%
  \renewcommand{\labelenumii}{(\theenumii)}%
}

\newcommand\Cinf{\mathcal{C}^\infty}

\def\ad{~\mbox{ad}~}
\def\char{~\mbox{char}~}
\def\Der{~\mbox{Der}~}
\def\Inder{~\mbox{Inder}~}
\def\End{~\mbox{End}~}
\def\Hom{~\mbox{Hom}~}
\def\Lie{~\mbox{Lie}~}

\begin{document} 
\maketitle 
\date
\begin{abstract} 
We classify simple linearly compact $n$-Lie superalgebras with $n>2$
over a field $\F$ of characteristic 0. The classification is based on a  
bijective correspondence between non-abelian
$n$-Lie superalgebras and transitive $\Z$-graded Lie
superalgebras of the form $L=\oplus_{j=-1}^{n-1} L_j$, 
where $\dim L_{n-1}=1$, $L_{-1}$ and $L_{n-1}$ generate $L$, and 
$[L_j, L_{n-j-1}] =0$ for all $j$, thereby reducing it 
to the known classification of simple linearly compact Lie superalgebras 
and their
$\Z$-gradings. The list consists of four examples, one of them being the
$n+1$-dimensional vector product $n$-Lie algebra, and the remaining three
infinite-dimensional $n$-Lie algebras.
\end{abstract} 
\section*{Introduction}


Given an integer $n \geq 2$, an $n$-Lie algebra $\fg$ is a
vector space over a field~$\FF$, endowed with an $n$-ary
anti-commutative product
\begin{displaymath}
  \Lambda^n \fg \to \fg \, , \quad a_1 
     \wedge \ldots \wedge a_n \mapsto [a_1,\ldots ,a_n]\, ,
\end{displaymath}
subject to the following Filippov-Jacobi identity:
\begin{equation}
\label{eq:0.1}
\begin{array}{c}
[a_1,\ldots ,a_{n-1}, [b_1,\ldots ,b_n]]= 
  [[a_1 ,\ldots ,a_{n-1}, b_1] , b_2 ,\ldots , b_n]
  + [b_1,[a_1,\ldots ,a_{n-1}, b_2], b_3 ,\ldots ,b_n]\\
+ \ldots
+ [b_1, \ldots ,b_{n-1}, [a_1 ,\ldots , a_{n-1},b_n]].
\end{array}
\end{equation}
The meaning of this identity is similar to that of the usual
Jacobi identity for a Lie algebra (which is a $2$-Lie algebra),
namely, given $a_1 , \ldots ,a_{n-1} \in \fg$, the map $D_{a_1,
  \ldots ,a_{n-1}} : \fg \mapsto \fg$, given by 
$D_{a_1,\ldots a_{n-1}}  (a) = [a_1,\ldots ,a_{n-1},a]$, 
is a derivation of the $n$-ary bracket.
These derivations are called inner.

The notion of an $n$-Lie algebra was introduced by V.T.\ Filippov  
in 1985 \cite{F}.
In this and several subsequent papers, \cite{F2}, \cite{Kas1}, \cite{Kas2}, \cite{L}, 
a structure theory of finite-dimensional $n$-Lie algebras over
a field $\FF$ of characteristic 0 was developed.  In
particular, W.X. Ling in \cite{L} discovered the following disappointing 
feature of $n$-Lie algebras for $n>2$:  there exists only
one simple finite-dimensional $n$-Lie algebra over an
algebraically closed field~$\FF$ of characteristic 0. It is given by
the vector product of $n$~vectors in the $n+1$-dimensional vector
space $V$, endowed with a non-degenerate symmetric bilinear form 
$(\cdot,\cdot)$.
Recall that, choosing dual bases $\{ a_i \}$ and $\{ a^i \}$ of~$V$,
i.e.,~$(a_i,a^j)=\delta_{ij}, i,j=1,...,n+1$, the vector product of
$n$~vectors from the basis $\{ a_i\}$ is defined as the following $n$-ary 
bracket:
\begin{displaymath}
  [a_{i_1},\ldots ,a_{i_n}]= \epsilon_{i_1, \ldots ,i_{n+1}}
  a^{i_{n+1}}\, ,
\end{displaymath}
where $\epsilon_{i_i ,\ldots ,i_{n+1}}$ is a non-zero totally
antisymmetric tensor with values in~$\FF$, and extended by n-linearity.  
This is a simple $n$-Lie algebra, which is called the vector product
$n$-Lie algebra; we denote it by $O^n$.  

Another example of an $n$-Lie
algebra appeared earlier in Nambu's generalization of
Hamiltonian dynamics \cite{N}.  
It is the space $\Cinf (M)$ of $\Cinf$-functions on a
finite-dimensional manifold~$M$, endowed with the following
$n$-ary bracket, associated to $n$~commuting vector fields
$D_1,\ldots,D_n$ on $M$:
\begin{equation}
  \label{eq:0.2}
  [f_1,\ldots ,f_n] = \det \left(
 \begin{array}{ccc}      
     D_1 (f_1) & \ldots & D_1 (f_n)\\       
\hdotsfor{3}\\ 
D_{n} (f_1) & \ldots & D_{n}(f_n)     
\end{array} \right) \, .
\end{equation}
The fact that this $n$-ary bracket satisfies the Filippov-Jacobi
identity was noticed  later by Filippov (who was unaware of
Nambu's work), and by Takhtajan \cite{T}, who introduced the
notion of an $n$-Poisson algebra (and was unaware of Filippov's
work).

A more recent important example of an $n$-Lie algebra structure
on $\Cinf (M)$, given by Dzhumadildaev \cite{D1}, is associated
to $n-1$ commuting vector fields $D_1 , \ldots ,D_{n-1}$ on $M$:
\begin{equation}
  \label{eq:0.3}
  [f_1 ,\ldots ,f_n] = \det \left( 
    \begin{array}{ccc}      
f_1 & \ldots & f_n\\      
D_1 (f_1) & \ldots & D_1 (f_n)\\       
\hdotsfor{3}\\ 
D_{n-1} (f_1) & \ldots & D_{n-1}(f_n)     
\end{array} \right) \, .
\end{equation}
In fact, Dzhumadildaev considered examples 
(\ref{eq:0.2}) and (\ref{eq:0.3}) 
in a
more general context, where $\Cinf (M)$ is replaced by an
arbitrary commutative associative algebra~$A$ over~$\FF$ and the
$D_i$ by derivations of $A$.  He showed in \cite{D4} that 
(\ref{eq:0.2}) and (\ref{eq:0.3}) 
satisfy the Filippov-Jacobi   identity if and only if
the vector space $\sum_i \FF D_i$ is closed under the Lie
bracket.

In the past few years there has been some interest in $n$-Lie
algebras in the physics community, related to $M$-branes in
string theory.  We shall quote here two sources --- a survey paper
\cite{FO}, containing a rather extensive list of references, and
a paper by Friedmann \cite{Fr}, where simple finite-dimensional
$3$-Lie algebras over $\CC$ were classified (she was unaware of
the earlier work).  At the same time we (the authors of the
present paper) have been completing our work \cite{CK4} on
simple rigid linearly compact superalgebras,  and it occurred to
us that the method of this work also applies to the
classification of simple linearly compact $n$-Lie superalgebras!

Our main result can be stated as follows.

\begin{theorem}
  \label{th:0.1}
\alphaparenlist
  \begin{enumerate}
  \item Any simple linearly compact $n$-Lie
    algebra with $n>2$, over an algebraically closed
    field~$\FF$ of characteristic 0, is isomorphic to one of the
following four examples:

\romanparenlistii
    \begin{enumerate}
    \item 
the $n+1$-dimensional vector product $n$-Lie algebra $O^n$;

\item 
the $n$-Lie algebra, denoted by $S^n$, which is the linearly compact
vector space of formal power
series $\FF [[x_1, \ldots ,x_n]]$, endowed with the $n$-ary bracket
(\ref{eq:0.2}), where $D_i = \frac{\partial}{\partial x_i}$;

\item 
the $n$-Lie algebra, denoted by $W^n$, which is the linearly compact vector
space of formal power
series $\FF [[x_1,\ldots ,x_{n-1}]]$, endowed with the
$n$-ary bracket (\ref{eq:0.3}), where $D_i = \frac{\partial}{\partial x_i}$;

\item 
  the $n$-Lie algebra, denoted by $SW^n$, which is the  direct
  sum of $n-1$ copies of $\FF [[x]]$, endowed with the following
  $n$-ary bracket, where $f^{\langle j \rangle}$ is an element of the
  $j$\st{th} copy and $f' = \frac{df}{dx}$:
\begin{align*}
  [f^{\langle j_1 \rangle}_1 , \ldots f^{\langle j_n \rangle}_n ] =0,
\hbox{\,\, unless \,\,} \{ j_1,\ldots ,j_n \}\supset \{ 1 ,\ldots , n-1 \},\\
  [f^{\langle 1 \rangle}_1 ,\ldots , f^{\langle k-1 \rangle}_{k-1},
    f^{\langle k \rangle}_k, f^{\langle k \rangle}_{k+1}, 
    f^{\langle k+1 \rangle}_{k+2}, \ldots , 
    f^{\langle n-1 \rangle}_n ]\\
    = (-1)^{k+n} (f_1 \ldots f_{k-1} (f'_kf_{k+1} - f'_{k+1}f_k)
       f_{k+2}\ldots f_n)^{\langle k \rangle}\, .
\end{align*}
    \end{enumerate}
\item 
There are no simple linearly compact $n$-Lie  superalgebras over
$\FF$, which are not $n$-Lie algebras, if $n>2$.

\end{enumerate}

\end{theorem}

Recall that a linearly compact algebra is a topological algebra,
whose underlying vector space is linearly compact, namely is a
topological product of finite-dimensional vector spaces, endowed
with discrete topology (and it is assumed that the algebra product is
continuous in this topology).  In particular, any
finite-dimensional algebra is automatically linearly compact.
The basic example of an infinite-dimensional linearly compact
space is the space of formal power series $\FF [[x_1 ,\ldots
,x_k]]$, endowed with the formal topology, or a direct sum of
a finite number of such spaces.

The proof of Theorem \ref{th:0.1} is based on a construction, which
associates to an $n$-Lie (super)algebra $\fg$ a pair $(\Lie
\fg$, $\mu$), where $\Lie\fg = \prod_{j \geq -1} \Lie_j\fg$
is a $\ZZ$-graded Lie superalgebra of depth 1 and  $\mu \in \Lie_{n-1}\fg$,
such that the following
properties hold:

\begin{list}{}{} 
\item 
(L1)\quad $\Lie\fg$ is transitive, i.e.,~if $a \in \Lie_j \fg$
with $j \geq 0$ and $[a,\Lie_{-1}\fg]=0$, then $a=0$;

\item 
(L2)\quad $\Lie \fg$ is generated by $\Lie_{-1}\fg$ and $\mu$;

\item 
  (L3)\quad  $[\mu ,\Lie_0\fg] =0$.

\end{list}

A pair $(L,\mu)$, where $L=\prod_{j \geq -1}L_j$ is a
transitive $\ZZ$-graded Lie superalgebra and $\mu \in L_{n-1}$,
such that (L2) and (L3) hold, is called
{\em admissible}.

The construction of the admissible pair $(\Lie \fg , \mu)$,
associated to an $n$-Lie (super)algebra $\fg$, uses the universal
$\ZZ$-graded Lie superalgebra $W(V) = \prod_{j \geq -1} W_j
(V)$, associated to a vector superspace $V$ (see Section \ref{Sec1}
for details).  One has
$W_j(V) = \Hom (S^{j+1}V,V)$, so that an element $\mu \in
W_{n-1} (V)$ defines a commutative $n$-superalgebra structure on $V$ 
and vice versa. {\em Universality} means that any transitive $\ZZ$-graded Lie
superalgebra $L= \prod_{j \geq -1} L_j$ with $L_{-1}=V$
canonically embeds in $W(V)$ (the embedding being given by $L_j
\ni a \mapsto \varphi_a \in W_j(V)$, where $\varphi_a (a_1,
\ldots ,a_{j+1}) = [\ldots [a,a_1],\ldots ,a_{j+1}]$). 

So, given a commutative $n$-ary product on a superspace $V$, we
get an element $\mu \in W_{n-1}(V)$, and we denote by $\Lie V$
the $\Z$-graded subalgebra of $W (V)$, generated by $W_{-1}(V)$ and
$\mu$.  The pair $(\Lie V,\mu)$ obviously satisfies properties
(L1) and (L2).

How do we pass from commutative to anti-commutative
$n$-superalgebras?  Given a commutative $n$-superalgebra $V$ with 
$n$-ary product $(a_1,\ldots ,a_n)$, the vector superspace $\Pi
V$ ($\Pi$ stands, as usual, for reversing the parity) becomes
an anti-commutative $n$-superalgebra with $n$-ary product
\begin{equation}
\label{eq:0.4}
[a_1 ,\ldots ,a_n] =p (a_1,\ldots ,a_n)(a_1,\ldots ,a_n),\\
\end{equation}
where 
\begin{equation}
\label{eq:0.5}
p(a_1,\ldots ,a_n)= \left\{
  \begin{array}{l}
    (-1)^{p(a_1)+ p(a_3)+ \cdots +p(a_{n-1})} \hbox{\,\, if \,\,}
    n \hbox{\,\,is even}\, \\
    (-1)^{p(a_2) +p(a_4)+ \cdots + p (a_{n-1})} \hbox{\,\, if \,\,}
      n \hbox{\,\,is odd}\, ,
  \end{array} \right.
\end{equation}
and vice versa.

Thus, given an anti-commutative $n$-superalgebra $\fg$, with
$n$-ary product $[a_1,\ldots ,a_n]$, we consider the vector
superspace $\Pi \fg$ with commutative $n$-ary product
$(a_1,\ldots ,a_n)$, given by (\ref{eq:0.4}), consider the element $\mu
\in W_{n-1} (\Pi \fg)$, 
 corresponding
to the latter $n$-ary product, and let $\Lie\fg$ be the graded
subalgebra of $W (\Pi \fg)$, generated by $W_{-1} (\Pi \fg)$
and $\mu$.

Note that properties (L1) and (L2) of the pair $(\Lie\fg,\mu)$
still hold, and it remains to note that property (L3) is
equivalent to the (super analogue of the) Filippov-Jacobi
identity.
Finally, the simplicity  of the $n$-Lie (super)algebra $\fg$ is
equivalent to 

\begin{list}{}{} 
\item 
(L4)\quad \hbox{the}\, $\Lie_0\fg$-module $\Lie_{-1}\fg$ 
\,\,\hbox{is irreducible}.
\end{list}

An admissible pair, satisfying property (L4) is called
irreducible.  Thus, the proof of Theorem~0.1 reduces to the
classification  of all irreducible admissible pairs $(L,\mu)$,
where $L$ is a linearly compact Lie superalgebra. It is not difficult to show,
as in \cite{CK4}, that $S\subseteq L \subseteq \Der S$, where 
$S$ is a simple linearly compact
Lie superalgebra and $\Der S$ is the Lie superalgebra of its
continuous derivations (it is at this point that the condition
$n>2$ is essential).

Up to now the arguments worked over an arbitrary field $\FF$.  In
the case $\FF$ is algebraically closed of characteristic $0$,
there is a complete classification of  simple linearly compact
Lie superalgebras, their derivations and their $\Z$-gradings
\cite{K1}, \cite{K2}, \cite{CK1}, \cite{CK2}.  Applying these
classifications completes the proof of Theorem~0.1

For example, if $\fg$ is a finite-dimensional simple $n$-Lie
superalgebra, then $\dim \Lie \fg < \infty$, and from \cite{K1} we see
that the only possibility for $\Lie \fg$ is $L=P/\FF1$, where
$P$ is the Lie superalgebra defined by the super Poisson bracket
on the Grassmann algebra in the indeterminates $\xi_1,\ldots
,\xi_{n+1}$, given by
\begin{equation}
\label{eq:0.6}
  \{ \xi_i, \xi_j \} = b_{ij}\, ,\,\,i,j=1,\ldots , n+1.
\end{equation}
where $(b_{ij})$ is a non-degenerate symmetric matrix, the
$\ZZ$-grading on $L$ being given by $\deg (\xi_{i_1}\ldots \xi_{i_s})=s-2$,
and $\mu = \xi_1 \xi_2 \ldots \xi_{n+1}$.  We conclude that $\fg$
is the vector product $n$-Lie algebra.
(The proof of this result in the non-super case, obtained by Ling \cite{L}, 
is based
on the study of the linear Lie algebra spanned by the derivations
$D_{a_1,...,a_{n-1}}$, and is applicable neither in the super nor
in the infinite-dimensional case.) 
We have no {\it a priori} proof of part~(b) of Theorem~0.1
 --- it comes out only after the classification process.

If $\char \FF =0$, we have a more precise result on the structure
of an admissible pair $(L,\mu)$.

\begin{theorem}
  \label{th:0.2}
If $\char \FF =0$ and  $(L,\mu)$ is an admissible pair, then $L =
\oplus^{n-1}_{j=-1} L_j$, where $L_{n-1}=\FF \mu\, , \,
S:=\oplus^{n-2}_{j=-1} L_j$ is an ideal in $L$, and $L_j =
(\ad L_{-1})^{n-j-1} \mu$ for $j=0, \ldots , n-1$.
\end{theorem}

Of course, Theorem \ref{th:0.2} reduces significantly the case wise
inspection in the proof of Theorem \ref{th:0.1}.  Moreover,  Theorem 
\ref{th:0.2}
can also be used in the study of representations of $n$-Lie
algebras.  Namely a representation of an $n$-Lie algebra
$\fg$ in a vector space~$M$ corresponds to an $n$-Lie
algebra structure on the semidirect product $L_{-1}=\fg \ltimes M$, 
where $M$ is an
abelian ideal.  Hence, by Theorem \ref{th:0.2}, we obtain a graded
representation of the Lie superalgebra $S = \oplus^{n-2}_{i=-1}
L_i$ in the graded supervector space $L(M) =
\sum^{n-2}_{i=-1}M_j$, $M_j = (\ad M)^{n-j-1} \mu$ so that $L_i
M_j \subset M_{i+j}$.  Such ``degenerate'' representations of the
Lie superalgebra $S$ are not difficult to classify, and this
corresponds to a classification of representations of the $n$-Lie
algebra $\fg$.  In particular, representations of the $n$-Lie algebra $O^n$
correspond to ``degenerate'' representations of the simple
Lie superalgebra $H(0,n)$ (finite-dimensional representations of $O^n$ were
classified in \cite{D3}, using Ling's method, mentioned above).  

Finally, note that, using our discussion on $\FF$-forms of simple
linearly compact Lie superalgebras in \cite{CK2}, we can extend
Theorem \ref{th:0.1} to the case of an arbitrary field $\FF$ of
characteristic~$0$.  The result is almost the same, namely the
$\FF$-forms are as follows: $O^n$, depending on the equivalence class of
the symmetric bilinear form up to a non-zero factor, and the $n$-Lie
algebras $S^n$, $W^n$ and $SW^n$ over $\FF$.

We would like to thank T. Friedmann and J. Thierry-Mieg, who drew our 
attention 
to this topic, and D.\ Balibanu, A.\ Dzhumadildaev, A.\ Kiselev and E.\ Zelmanov for useful 
discussions and correspondence.
Also we wish to thank the referees for very useful comments.
In particular, one of the referees pointed out that in the case
$n=3$ the construction of \cite{MFM} is closely related to our construction.

\section{Preliminaries on $n$-superalgebras}\label{Sec1}
\label{sec:1}

Let $V$ be a vector superspace over a field $\FF$, namely we
have a decomposition $V =V_{\bar{0}} \oplus V_{\bar{1}}$ in a
direct sum of subspaces, where $V_{\bar{0}}$
(resp.\ $V_{\bar{1}}$) is called the subspace of even (resp.\ odd)
elements; if $v \in V_\alpha$, $\alpha \in \ZZ /2\ZZ = \{
\bar{0}, \bar{1} \}$ we write $p (v) =\alpha$.  Given two
vector superspaces $U$ and $V$, the space $\Hom (U,V)$ is naturally
a vector superspace, for which even (resp.\ odd) elements are parity
preserving (resp.\ reversing) maps; also $U \otimes V$ is a vector
superspace via letting $p(a \otimes b) = p(a) + p(b)$ for $a \in
U$, $b \in V$.

In particular, the tensor algebra $T(V) = \oplus_{j \in \ZZ_+}
T^j(V)$ is  an associative superalgebra.  The symmetric
(resp.\ exterior) superalgebra over $V$ is the quotient of the
superalgebra $T(V)$ by the $2$-sided ideal, generated by the
elements $a \otimes b - (-1)^{p(a)p(b)} b \otimes a$ (resp.\ $a
\otimes b + (-1)^{p(a)p(b)} b \otimes a$), where $a,b \in V$.
They are denoted by $S(V)$ and $\Lambda (V)$ respectively.  Both
inherit a $\ZZ$-grading from $T(V)$:  $S(V) = \oplus_{j \in \ZZ_+}
S^j$, $\Lambda (V) = \oplus_{j \in \ZZ_+} \Lambda^j(V)$.

A well known trivial, but important, observation is that the
reversal of parity of $V$, i.e.,~taking the superspace $\Pi V$,
where $(\Pi V )_\alpha = V_{\alpha + \bar{1}}$, establishes a
canonical isomorphism:
\begin{equation}
  \label{eq:1.1}
  S (\Pi V) \simeq \Lambda (V)\, .
\end{equation}

\begin{definition}
  \label{def:1.1}

Let $n \in \ZZ_+$ and let $V$ be a vector superspace.  An
$n$-superalgebra structure (or $n$-ary product) on $V$ of parity
$\alpha \in \ZZ /2\ZZ$ is a linear map $\mu : T^n (V) \to V$ of
parity $\alpha$.  A commutative (resp.\ anti-commutative)
$n$-superalgebra of parity $\alpha$ is a linear map $\mu : S^n
(V) \to V$ (resp.\ $\Lambda^n V \to V$) of parity $\alpha$,
denoted by $\mu (a_1 \otimes \cdots \otimes a_n) = (a_1,\ldots
,a_n)$ (resp.\ $= [a_1, \ldots ,a_n]$).

\end{definition}

\begin{lemma}
  \label{lem:1.2}

Let $(V,\mu)$ be an anti-commutative $n$-superalgebra.  Then
$(\Pi V ,\bar{\mu})$ is a commutative $n$-superalgebra (of
parity $p(\mu)+n-1 ~mod~ 2$) with the $n$-ary product
\begin{equation}
  \label{eq:1.2}
  \bar{\mu} (a_1 \otimes \cdots \otimes a_n) = p(a_1,\ldots,a_n)
    \mu(a_1\otimes \ldots \otimes a_n),
\end{equation}
where
\begin{equation}
\label{eq:1.3}
p(a_1,\ldots ,a_n) =
(-1)^{\sum_{k=0}^{[\frac{n-2}{2}]}p(a_{n-1-2k})},
\end{equation}
and vice versa.
\end{lemma}

{\bf Proof.} Denote by $p'$ the parity in $\Pi V$. Then,
for $a_1,\dots, a_n\in V$,
$p'(\bar{\mu}(a_1\otimes\dots\otimes a_n))=p(\mu)+1+\sum_{i=1}^np(a_i)=
\sum_{i=1}^np'(a_i)+p(\mu)+n+1 ~~\mbox{mod} ~2$, i.e., 
$p'(\bar{\mu})=p(\mu)+n-1 ~~\mbox{mod} ~2$.
Besides, we have $p(a_1,\dots,a_i,a_{i+1},\dots,a_n)
p(a_1,\dots,a_{i+1},a_i,\dots,a_n)=(-1)^{p(a_i)+p(a_{i+1})}$, hence:
$$\bar{\mu}(a_1\otimes\dots\otimes a_i\otimes a_{i+1}\otimes \dots\otimes a_n)
=p(a_1,\dots,a_n)\mu(a_1\otimes\dots\otimes a_i\otimes a_{i+1}\otimes\dots
\otimes a_n)$$
$$=-(-1)^{p(a_i)p(a_{i+1})}p(a_1,\dots,a_n)\mu(a_1\otimes\dots\otimes 
a_{i+1}\otimes a_i\otimes \dots\otimes a_n)$$
$$=-(-1)^{p(a_i)p(a_{i+1})}p(a_1,\dots,a_n)p(a_1,\dots,a_{i+1}, a_i,
\dots, a_n)\bar{\mu}(a_1\otimes\dots\otimes 
a_{i+1}\otimes a_i\otimes \dots\otimes a_n)$$
$$=(-1)^{p'(a_i)p'(a_{i+1})}\bar{\mu}(a_1\otimes\dots\otimes 
a_{i+1}\otimes a_i\otimes \dots\otimes a_n).$$
\hfill$\Box$

\begin{definition}
  \label{def:1.2}
A derivation $D$ of parity $\alpha \in \ZZ /2\ZZ$ of an
$n$-superalgebra $(V,\mu)$ is an endomorphism of the vector
superspace $V$ of parity $\alpha$, such that:
$$  D (\mu (a_1 \otimes \cdots \otimes a_n)) = (-1)^{\alpha  p(\mu)}
    (\mu (D a_1 \otimes a_2 \otimes \cdots \otimes a_n)
    + (-1)^{\alpha p (a_1)} \mu (a_1 \otimes Da_2 \otimes \cdots
    \otimes a_n) + \cdots $$
%
$$+ (-1)^{\alpha (p(a_1) + \cdots + p (a_{n-1}))}
      \mu (a_1 \otimes \cdots \otimes D (a_n)))\, .$$
%

\end{definition}

It is clear that derivations of an $n$-superalgebra $V$ form a Lie 
superalgebra, which is denoted by $\Der V$. It is not difficult to show that
all inner derivations of an $n$-Lie algebra $\fg$ span an ideal of $\Der \fg$
(see e.g. \cite{D3}), denoted by $\Inder \fg$. 

Now we recall the construction of the universal Lie superalgebra $W(V)$,
associated to the vector superspace $V$.  For an integer $k \geq
-1$ let $W_k (V) = \Hom (S^{k+1}(V),V)$, in other words ,
$W_k (V)$ is the vector superspace of all commutative
$k+1$-superalgebra structures on $V$, in
particular, $W_{-1} (V) =V$, $W_0 (V) = \End (V) $, $W_1 (V)$ is
the space of all commutative superalgebra structures on $V$,
etc.  We endow the vector superspace
\begin{displaymath}
  W (V) = \prod^\infty_{k=-1} W_k (V)
\end{displaymath}
with a product $f \Box g$, making $W(V)$ a $\ZZ$-graded
superalgebra, given by the following formula for $f \in W_p (V)$,
$g \in W_q(V)$:
\begin{equation}
  \label{eq:1.5}
f\Box g(x_0, \dots, x_{p+q})=
\sum_{\substack
{i_0<\dots <i_q\\
i_{q+1}<\dots< i_{p+q}}}\epsilon(i_0,\dots, i_q, i_{q+1},\dots,
i_{p+q})f(g(x_{i_0},\dots, x_{i_q}), x_{i_{q+1}},\dots, x_{i_{p+q}}),
\end{equation}
where $\epsilon = (-1)^N$, $N$ being the number of interchanges
of indices of odd $x_i$'s in the permutation $\sigma (s) =i_s$,
$s=0,1,\ldots ,p+q$.  Then the bracket 
\begin{equation}
  \label{EQ:1.5}
  [f,g] = f \Box g - (-1)^{p(f)p(q)} g \Box f
\end{equation}
defines a Lie superalgebra structure on $W(V)$.

\begin{lemma}
  \label{lem:1.4}

Let $V$  be a vector superspace and let $\mu \in W_{n-1} (V)$,
$D \in W_0 (V)$.  Then

\begin{list}{}{}
\item (a) \quad $[\mu ,D]=0$ if and only if $D$ is a derivation
  of the $n$-superalgebra $(V,\mu)$.

\item (b) \quad 
$D$ is a derivation of parity $\alpha$ of the commutative $n$-superalgebra 
$(V ,{\mu})$ if and only
if $D$ is a derivation of parity $\alpha$ of the anti-commutative $n$-
superalgebra $(\Pi V ,\bar{\mu})$, where 
  \begin{displaymath}
    \bar{\mu} (a_1 \otimes \cdots \otimes a_n) 
       = p (a_1,\ldots ,a_n) \mu (a_1 \otimes \cdots \otimes a_n)\,,
  \end{displaymath}
and $p(a_1,\ldots ,a_n)$ is defined by (\ref{eq:1.3}).
\end{list}

\end{lemma}

{\bf Proof.} By (\ref{EQ:1.5}) and (\ref{eq:1.5}), we have:
$$[\mu,D](b_1\otimes\dots\otimes b_n)=(\mu\Box D)(b_1\otimes\dots\otimes b_n)-(-1)^{\alpha p(\mu)}
(D\Box\mu)(b_1\otimes\dots\otimes b_n)=$$
$$\sum_{\substack
{i_1\\
i_2<\dots<i_n}}(\varepsilon(i_1,\dots,i_n)\mu(D(b_{i_1})\otimes b_{i_2}\otimes
\dots\otimes b_{i_n}))-(-1)^{\alpha p(\mu)}D(\mu(b_1\otimes\dots\otimes b_n)),$$
where $\alpha$ is the parity of $D$.
Therefore $[\mu, D]=0$ if and only if
$$D(\mu(b_1\otimes\dots\otimes b_n))=(-1)^{\alpha p(\mu)}\sum_{\substack
{i_1\\
i_2<\dots<i_n}}(\varepsilon(i_1,\dots,i_n)\mu(D(b_{i_1})\otimes 
b_{i_2}\otimes
\dots\otimes b_{i_n}))$$
$$=(-1)^{\alpha p(\mu)}(\mu(D(b_1)\otimes b_2\otimes\dots\otimes b_{n-1})+
(-1)^{\alpha p(b_1)}\mu(b_1\otimes D(b_2)\otimes\dots\otimes b_{n-1})+
\dots$$
$$+(-1)^{\alpha(p(b_1)+\dots +p(b_{n-1}))}\mu(b_1\otimes\dots\otimes b_{n-1}\otimes D(b_n))$$
i.e., if and only if $D$ is a derivation of $(V,\mu)$ of parity $\alpha$,
proving $(a)$. 

In order to prove $(b)$, note that $D$ is a derivation
of parity $\alpha$ of $(V,\mu)$ if and only if
$$D(\bar{\mu}(a_1\otimes \dots\otimes a_n)=p(a_1,\dots, a_n)D(\mu(a_1\otimes \dots
\otimes a_n))
=p(a_1,\dots,a_n)((-1)^{\alpha p(\mu)}(\mu(D(a_1)\otimes a_2\otimes\dots
a_n)$$
$$+(-1)^{\alpha p(a_1)}\mu(a_1\otimes D(a_2)\otimes\dots a_n)+
\dots 
+(-1)^{\alpha(p(a_1)+\dots + p(a_{n-1}))}\mu(a_1\otimes \dots\otimes D(a_n)))$$
$$=p(a_1,\dots,a_n)(-1)^{\alpha p(\mu)}p(D(a_1),a_2\dots, a_n)
({\bar{\mu}}(D(a_1)\otimes a_2\dots \otimes a_n)$$
$$+(-1)^{\alpha p(a_1)}p(D(a_1),a_2,\dots, a_n)
p(a_1,D(a_2),\dots, a_n)
{\bar{\mu}}(a_1\otimes D(a_2)\dots \otimes a_n)+\dots$$
$$+(-1)^{\alpha (p(a_1)+\dots +p(a_{n-1}))}
p(D(a_1), a_2,\dots, a_n)p(a_1,\dots, D(a_n))
{\bar{\mu}}(a_1\otimes \dots \otimes D(a_n)))$$
If $n$ is even, we have:
$$\begin{array}{c}
p(a_1,\dots,a_n)(-1)^{\alpha p(\mu)}p(D(a_1), a_2\dots,a_n)=
(-1)^{\alpha(p(\mu)+1)},\\
(-1)^{\alpha p(a_1)}p(D(a_1),a_2\dots,a_n)
p(a_1,D(a_2),\dots, a_n)=(-1)^{\alpha(p(a_1)+1)},\\
\vdots\\
(-1)^{\alpha (p(a_1)+\dots +p(a_{n-1}))}p(D(a_1),a_2\dots,a_n)
p(a_1,\dots, D(a_n))=
(-1)^{\alpha (p(a_1)+\dots +p(a_{n-1})+1)}
\end{array}
$$
If $n$ is odd, we have:
$$\begin{array}{c}
p(a_1,\dots,a_n)(-1)^{\alpha p(\mu)}p(D(a_1),a_2\dots, a_n)=
(-1)^{\alpha p(\mu)},\\
(-1)^{\alpha p(a_1)}p(D(a_1),a_2\dots,a_n)
p(a_1,D(a_2),\dots, a_n)=(-1)^{\alpha(p(a_1)+1)},\\
\vdots\\
(-1)^{\alpha (p(a_1)+\dots +p(a_{n-1}))}
p(D(a_1),a_2\dots, a_n)p(a_1,\dots, D(a_n))
(-1)^{\alpha (p(a_1)+\dots +p(a_{n-1}))}.
\end{array}
$$
Since $\bar{\mu}$ has parity equal to $p(\mu)+n-1$ mod $2$, $(b)$ is proved.
\hfill$\Box$

\section{The main construction}
\label{sec:2}

Let $(\fg ,\mu)$ be an anti-commutative $n$-superalgebra over a field $\FF$ with
$n$-ary product $[a_1 ,\ldots , a_n]$, and let $V = \Pi \fg$.
Consider the universal Lie superalgebra $W (V) =
\prod^\infty_{k=-1} W_k (V)$, and let $\bar{\mu} \in W_{n-1}
(V)$ be the element defined by (\ref{eq:1.2}).  Let $\Lie \fg =
\prod^\infty_{j=-1} \Lie_j\fg$ be the $\ZZ$-graded subalgebra of
the Lie superalgebra $W(V)$, generated by $W_{-1}(V) = V$ and
$\bar{\mu}$.

\begin{lemma}
  \label{lem:2.1}
\alphaparenlist
  \begin{enumerate}
  \item 
      $\Lie \fg$ is a transitive subalgebra of $W(V)$.

\item 
  If $D \in \Lie_0\fg$, then the action of $D$ on $\Lie_{-1}\fg
  =V (= \Pi \fg)$ is a derivation of the $n$-superalgebra $\fg$ 
if and only if $[D,\bar{\mu}] =0$.

\item 
  $\Lie_0\fg$ is generated by elements of the form
  \begin{equation}
    \label{eq:2.1}
    (\ad a_1) \ldots (\ad a_{n-1}) \bar{\mu}\, , 
    \hbox{\,\,where \,\,} a_i\in\Lie_{-1} \fg =V \, .
  \end{equation}

  \end{enumerate}

\end{lemma}

{\bf Proof.} (a) is clear since $W(V)$ is transitive and the latter holds
since, for $f \in W_k (V)$ and $a,a_1 ,\ldots ,a_k \in W_{-1} (V)
=V$ one has:
\begin{displaymath}
  [f,a] (a_1,\ldots ,a_k) = f (a,a_1,\ldots ,a_k)\, .
\end{displaymath}
(b) follows from Lemma \ref{lem:1.4}.

In order to prove (c) let 
$\tilde{L}_{-1}=V$ and let $\tilde{L}_0$ 
be the subalgebra of the
Lie algebra $W_0 (V)$, generated by elements (2.1).  Let
$\prod_{j \geq -1} \tilde{L}_j$ be the full
prolongation of $\tilde{L}_{-1} \oplus \tilde{L}_0$,
i.e.,~$\tilde{L}_j = \{ a \in W_j (V) | [a,\tilde{L}_{-1}]
\subset \tilde{L}_{j-1}\}$ for $j \geq 1$.  This is a subalgebra
of $W(V)$, containing $V$ and $\bar{\mu}$, hence $\Lie \fg$.
This proves (c).
\hfill$\Box$

\begin{definition}
\label{def:2.2}
An $n$-Lie superalgebra is an anti-commutative $n$-superalgebra $\fg$ of 
parity $\alpha$, such that all endomorphisms $D_{a_1,\ldots ,a_{n-1}}$ of $\fg$
($a_1,\ldots a_{n-1} \in \fg$), defined by

\begin{displaymath}
D_{a_1,\ldots a_{n-1}}(a) = [a_1,\ldots ,a_{n-1},a],
\end{displaymath}
are derivations of $\fg$, i.e., the following Filippov-Jacobi identity holds:
\begin{equation}
\label{FJ}
\begin{array}{c}
[a_1, \dots , a_{n-1},[b_1, \dots , b_n]]=
(-1)^{\alpha (p(a_1)+\dots +p(a_{n-1}))}([[a_1, \dots , a_{n-1}, b_1], b_2,
\dots , b_n]+\\
(-1)^{p(b_1)(p(a_1)+\dots +p(a_{n-1}))}[b_1, [a_1, \dots , a_{n-1}, b_2], 
b_3,\dots , b_n]+\\
\dots +(-1)^{(p(b_1)+\dots +p(b_{n-1}))(p(a_1)+\dots +
p(a_{n-1}))}[b_1,\dots , b_{n-1}, [a_1, \dots , a_{n-1}, b_n]]).
\end{array}
\end{equation}
\end{definition}

Recall from the introduction that the pair $(L,\mu)$, where
$L=\prod_{j \geq -1}L_j$ is a $\ZZ$-graded Lie superalgebra and $\mu \in 
L_{n-1}$, 
is called admissible if 
properties (L1),(L2) and (L3) hold. Two admissible pairs 
$(L,\mu)$ and $(L',\mu')$ are called isomorphic if there exists a Lie
superalgebra isomorphism 
$\phi:L \mapsto L'$, 
such that $\phi:L_j= L'_j$ for all $j$ and $\phi(\mu) \in \FF^{\times}\mu'$.

The following corollary of Lemma 2.1 is immediate. 

\begin{corollary}
\label{cor:2.3}
If $\g$ is an $n$-Lie superalgebra, then the pair $(\Lie \fg,\bar{\mu})$
is admissible.
\end{corollary}

Now it is easy to prove the following key result.

\begin{proposition}
\label{prop:2.4}
The map $\fg \mapsto (\Lie \fg, \bar{\mu})$ induces a bijection between 
isomorphism classes of $n$-Lie superalgebras, considered up to rescaling the 
$n$-ary bracket, and isomorphism classes of admissible pairs. Under this 
bijection, simple $n$-Lie algebras correspond to irreducible admissible pairs.
Moreover,
$\fg$ is linearly compact if and only if $\Lie\fg$ is.
\end{proposition}
{\bf Proof.}
Given an admissible pair $(L,\mu)$, where $L=\prod_{j \geq -1}L_j$, $\mu \in
L_{n-1}$, we let $\fg=\Pi L_{-1}$, and define an $n$-ary bracket on $\fg$
by the formula
$$ [a_1 ,\ldots ,a_n] =p (a_1,\ldots ,a_n)[...[\mu,a_1]\ldots , a_n],$$
where $p(a_1, \ldots, a_n)$ is given by (\ref{eq:1.3}).
Obviously, this $n$-ary bracket is anti-commutative. The Filippov-Jacobi 
identity follows from the property (L3) using the embedding of $L$
in $W(L_{-1})$ and applying Lemma 2.1(b).
Thus, $\fg$ is an $n$-Lie superalgebra. Due to properties (L1) and (L2),
we obtain the bijection of the map in question. 
It is obvious that $\fg$ is simple if and only if the pair $(L,\mu)$ is irreducible.
The fact that the linear 
compactness of $\fg$ implies that of $\Lie \fg$ is proved in the same way
as Proposition 7.2(c) from \cite{CK4}.
\hfill$\Box$

\begin{remark}\em 
\label{rem:2.5}
If $\fg$ is a finite-dimensional $n$-Lie algebra, then $\Lie_{-1}\fg=\Pi\fg$
is purely odd, hence 
$\dim W(\Pi \fg)<\infty$ 
and therefore $\dim \Lie \fg<\infty$. In the super case
this follows from Theorem \ref{th:0.2} if $\char \F=0$,
and from the fact that any finite-dimensional subspace of $W(V)$
generates a finite-dimensional subalgebra if $\char \F >0$.
Thus, an $n$-Lie superalgebra $\fg$ is finite-dimensional
if and only if the Lie superalgebra $\Lie \fg$ is
finite-dimensional.
\end{remark}

\begin{remark}\em 
\label{rem:2.6}
Let $V$ be a vector superspace. Recall that a sequence of anti-commutative
$(n+1)$-ary products $d_n$ , $n=0,1,\ldots$, of parity $n+1 \mod 2$ on $V$
endow $V$ with a structure of a 
homotopy Lie algebra if they satisfy a sequence of certain quadratic 
identities, which mean that $d_0^2=0$, $d_1$ is a Lie (super)algebra
bracket modulo the image of $d_0$ and $d_0$ is the derivation of
this bracket, etc \cite{SS}. (Usually one also requires a $\Z$-grading
on $V$ for which $d_n$ has degree $n-1$, but we ignore this requirement 
here.) On the other hand, recall that if $\mu_n$ is an 
$(n+1)$-ary anti-commutative product on $V$ of parity $n+1\mod 2$,
then the $(n+1)$-ary product 
$\bar{\mu}_n$, defined in Lemma \ref{lem:1.4},
is a commutative odd product 
on the vector superspace $\Pi V$. It is easy to see that the sequence of
$(n+1)$-ary products $\bar{\mu}_n$ define a homotopy Lie algebra 
structure on $\Pi V$ if and only if the odd element $\mu=\sum_n \bar{\mu}_n
\in W(\Pi V)$ satisfies the identity $[\mu,\mu]=0$. As above, we can
associate to a given homotopy Lie algebra structure on $V$ the subalgebra
of $W(\Pi V)$, denoted by $\Lie (V,\mu)$, which is generated by 
$W_{-1}(\Pi V)$ and all $\bar{\mu}_n$, $n=0,1,\ldots$. If the superspace $V$
is linearly compact and the homotopy Lie algebra is simple with 
$\bar{\mu}_n \neq 0$ for some $n>2$,
then the derived algebra of $\Lie (V,\mu)$ is simple, hence
is of one of the types $X(m,n)$, according to the classification of \cite{K2}.
Then the simple homotopy Lie algebra is called of type $X(m,n)$.
(Of course, there are many homotopy Lie algebras of a given type.)
Lemma \ref{Askar} below shows, in particular, that in characteristic 0
any $n$-Lie superalgebra of parity $n\mod 2$ is a homotopy Lie algebra, 
for which $\bar{\mu}_j=0$
if $j \neq n-1$. This was proved earlier in \cite{D2} and \cite{MV}.
\end{remark} 

\section{Proof of Theorem \ref{th:0.2}}
\label{sec:3}
For the sake of simplicity we consider the $n$-Lie algebra case, i.e. we 
assume that
$L_{-1}$ is purely odd. The same proof works verbatim when $L_{-1}$ is not 
purely odd, using identity (\ref{FJ}). Alternatively,
the use of the standard Grassmann envelope argument reduces the 
case of $n$-Lie superalgebras of even parity to the case of $n$-Lie algebras.

First, introduce some notation. Let $S_{2n-1}$ be the group of permutations
of the $2n-1$ element set $\{1,\ldots,2n-1\}$ and, for $\sigma\in S_{2n-1}$, 
let $\varepsilon(\sigma)$ be the sign of $\sigma$. 
Denote by $S$ the subset
of $S_{2n-1}$ consisting of permutations $\sigma$, such that 
$\sigma(1)<\dots<\sigma(n-1)$,
$\sigma(n)<\dots<\sigma(2n-1)$. Consider the following subsets of $S$
($l$ and $s$ stand for ``long'' and ``short'' as in \cite{D2}):
$$S^{l_1}=\{\sigma\in S|~\sigma(2n-1)=2n-1\},$$
$$S^{s_1}=\{\sigma\in S|~\sigma(n-1)=2n-1\}.$$
It is immediate to see that $S=S^{l_1}\cup S^{s_1}$.
Likewise, let
$$S^{l_1l_2}=\{\sigma\in S^{l_1}~|~\sigma(2n-2)=2n-2\},$$
$$S^{l_1s_2}=\{\sigma\in S^{l_1}~|~\sigma(n-1)=2n-2\},$$
$$S^{s_1l_2}=\{\sigma\in S^{s_1}~|~\sigma(2n-1)=2n-2\},$$
$$S^{s_1s_2}=\{\sigma\in S^{s_1}~|~\sigma(n-2)=2n-2\}.$$
Then $S^{l_1}=S^{l_1l_2}\cup S^{l_1s_2}$, and
$S^{s_1}=S^{s_1l_2}\cup S^{s_1s_2}$.
Likewise, we define the subsets $S^{a_1\dots a_k}$, with $a=l$ or
$a=s$,
for $1\leq k\leq 2n-1$, so that
\begin{equation}
\label{l+s}
S^{a_1\dots a_{k-1}}=
S^{a_1\dots a_{k-1}s_k}\cup
S^{a_1\dots a_{k-1}l_k}.
\end{equation}

\begin{lemma} \label{Askar}
 If $(L=\prod_{j\geq -1}L_j,\mu)$ is an admissible pair,
then $[(ad L_{-1})^{n-j-1}\mu,\mu]=0$ for every $j=0,\dots,n-1$.
\end{lemma}
{\bf Proof.}
If $(L,\mu)$ is an admissible pair, then, by Lemma \ref{lem:2.1}$(c)$,
$L_0=(ad L_{-1})^{n-1}\mu$,  hence
$[(ad L_{-1})^{n-1}$ $\mu, \mu]=0$ by property (L3).
Now we will show that $[(ad L_{-1})^{n-j-1}\mu,\mu]=0$ for every 
$j=1,\dots,n-1$.  By Lemma \ref{lem:2.1}$(b)$ and property (L3),
the Filippov-Jacobi identity holds for elements in $\Pi L_{-1}$
with product (\ref{eq:1.2}).
Let $x_1, \dots, x_{n-j-1}\in L_{-1}$.
By definition,
$[\mu, [x_1,\dots, [x_{n-j-1},\mu]]]=
\mu\Box [x_1,\dots, [x_{n-j-1},\mu]]-(-1)^{j(n-1)}$
$[x_1,\dots, [x_{n-j-1},\mu]]\Box\mu.$
One checks by a direct calculation, using the Filippov-
Jacobi identity, that 
$$\mu\Box Alt [x_1,\dots, [x_{n-j-1},\mu]]-
(-1)^{j(n-1)}
[x_1,\dots, [x_{n-j-1},\mu]]\Box\mu=0$$
where, for $f\in Hom(V^{\otimes k},V)$, $Alt f\in Hom(V^{\otimes k},V)$ 
denotes the alternator
of $f$, i.e., $Alt f(a_1,\dots, a_k)$ $=\sum_{\sigma\in S_k}
f(a_{\sigma(1)},$ $\dots, a_{\sigma(k)})$. Hence, since $\mu\in W_{n-1}(L_{-1})$, we have
$$[\mu, [x_1,\dots, [x_{n-j-1},\mu]]]=
((j+1)!+1)\mu\Box [x_1,\dots, [x_{n-j-1},\mu]].$$
We will prove a stronger statement than the lemma, namely, we will
show that, for every $j=1,\dots,n-1$ one has: 
\begin{equation}
\label{EQ:3.2}
\mu\Box [x_1,\dots, [x_{n-j-1},\mu]]=0.
\end{equation}
Note that, by definition, for $a_1,\dots, a_{n+j}\in L_{-1}$, we have:
$$(\mu\Box [x_1,\dots, [x_{n-j-1},\mu]])(a_1,\dots, a_{n+j})=
\sum_{\substack
{\sigma(1)<\dots <\sigma(j+1)\\
\sigma(j+2)<\dots<\sigma(n+j)
}}\varepsilon(\sigma)\mu([x_1,\dots,[x_{n-j-1},\mu]]
(a_{\sigma(1)},\dots, a_{\sigma(j+1)}),$$
$$ a_{\sigma(j+2)},\dots,
a_{\sigma(n+j)})=
\sum_{\substack
{\sigma(1)<\dots <\sigma(j+1)\\
\sigma(j+2)<\dots<\sigma(n+j)
}}\varepsilon(\sigma)\mu(\mu(x_{n-j-1},\dots, x_1,
a_{\sigma(1)},\dots, a_{\sigma(j+1)}),a_{\sigma(j+2)},\dots, a_{\sigma(n+j)}).
$$
Therefore (\ref{EQ:3.2}) is equivalent to the following:
\begin{equation}
\sum_{
\sigma\in S^{l_1,\dots,l_{n-j-1}}
}\varepsilon(\sigma)\mu(x_{\sigma(1)},\dots, x_{\sigma(n-1)},\mu(x_{\sigma(n)},\dots,x_{\sigma(2n-1)}))=0.
\label{EQ:3.1}
\end{equation}

Set
$A_{\sigma}=
\varepsilon(\sigma)\mu(x_{\sigma(1)},\dots, x_{\sigma(n-1)},\mu(x_{\sigma(n)},\dots,x_{\sigma(2n-1)}))$, 
$Q_{l_1\dots l_t}=
\sum_{\sigma\in S^{l_1\dots l_t}}A_{\sigma}$, and similarly define
$Q_{a_1\dots a_t}$, where $a=s$ or $a=l$. Then (\ref{EQ:3.1}) is equivalent to
$Q_{l_1\dots l_{n-j-1}}=0$. In fact, we shall prove more:
\begin{equation}
\label{thesis}
Q_{l_1\dots l_t}=Q_{l_1\dots l_{t-1}s_t}=Q_{s_1\dots s_{t-1}l_t}
=Q_{s_1\dots s_t}=0 ~~~\mbox{for}~ t=0,\dots, n-2.
\end{equation}
For $t=0$ and $t=1$, equality (\ref{thesis}) can be proved as
in \cite[Proposition 2.1]{D2}. Namely, by the Filippov-Jacobi identity,
for any $\sigma\in S^{s_1}$, $A_{\sigma}$ can be written as a
sum of $n$ elements $A_{\tau}$, where $\tau\in S^{l_1}$ is such that
$\{\tau(1),\dots, \tau(n-1)\}\subset\{\sigma(n),\dots, \sigma(2n-1)\}$.
Since the sets $\{\tau(1),\dots, \tau(n-1)\}$
and $\{\sigma(n),\dots, \sigma(2n-1)\}$ have $n-1$ and $n$ elements,
respectively, there exists only one $i$ such that
$\{\tau(1),\dots, \tau(n-1)\}\cup \{i\}=\{\sigma(n),\dots, \sigma(2n-1)\}$.
Then $i\leq 2n-2$ and $i\neq \tau(1),\dots, \tau(n-1)$. Therefore
there are $n-1$ possibilities to choose $i$. It follows that
\begin{equation}
\label{Askar1}
Q_{s_1}=(n-1)Q_{l_1}.
\end{equation}
Likewise, by the Filippov-Jacobi identity, for any $\sigma\in S^{l_1}$, $A_{\sigma}$ can be written as a
sum of one element $A_{\rho}$ with $\rho\in S^{l_1}$,
and $n-1$ elements $A_{\tau}$, with $\tau\in S^{s_1}$ such that
$\{\tau(1),\dots, \tau(n-2)\}\subset\{\sigma(n),\dots, \sigma(2n-2)\}$.
As above, there exists only one $i$ such that
$\{\tau(1),\dots, \tau(n-2)\}\cup \{i\}=\{\sigma(n),\dots, \sigma(2n-2)\}$.
Such an $i$ is different from $\tau(1),\dots, \tau(n-2)$
and $i\leq 2n-2$. Hence
there are $n$ possibilities to choose $i$. Notice that
\begin{equation}
\label{rho}
\rho=\left(\begin{array}{ccccccc}
1 & \dots & n-1 & n & \dots & 2n-2 & 2n-1\\
\sigma(n) & \dots & \sigma(2n-2) & \sigma(1) & \dots & \sigma(n-1) & 2n-1
\end{array}
\right),
\end{equation}
hence $\varepsilon(\rho)=(-1)^{n-1}\varepsilon({\sigma})$. Therefore
\begin{equation}
\label{Askar2}
Q_{l_1}=nQ_{s_1}+(-1)^{n-1}Q_{l_1}.
\end{equation}
Equations (\ref{Askar1}) and (\ref{Askar2}) form a system of
two linear equations in the two indeterminates $Q_{s_1}$ and $Q_{l_1}$,
whose determinant is equal to $n^2-n-1-(-1)^n$, which is
 different from zero for every $n>2$. It follows that $Q_{s_1}=0=Q_{l_1}$, i.e., 
(\ref{thesis}) is proved for $t=1$. Since, as we have already noticed,
$S=S^{l_1}+S^{s_1}$, (\ref{thesis}) 
 for $t=0$ also follows.

Now we argue by induction on $t$. We already proved (\ref{thesis})
for $t=0$ and $t=1$. Assume that
$$Q_{l_1\dots l_{t-1}}=Q_{l_1\dots l_{t-2}s_{t-1}}=
Q_{s_1\dots s_{t-2}l_{t-1}}=Q_{s_1\dots s_{t-1}}=0$$
for some $1\leq t<n-2$.
Similarly as above, by the Filippov-Jacobi identity, for any $\sigma\in S^{s_1\dots s_t}$, 
$A_\sigma$ can be
written as a sum of $n$ elements $A_{\tau}$ with
$\tau\in S^{l_1\dots l_t}$, such that
$\{\tau(1),\dots,\tau(n-1)\}\subset \{\sigma(n),\dots, \sigma(2n-1)\}$,
i.e.,
$\{\tau(1),\dots,\tau(n-1)\}\cup \{i\}=
 \{\sigma(n),\dots, \sigma(2n-1)\}$, for some $i\leq 2n-t-1$ 
and $i\neq \tau(1),
\dots, \tau(n-1)$. It follows that there are $n-t$ choices for $i$, hence
\begin{equation}
Q_{s_1\dots s_t}=(n-t)
Q_{l_1\dots l_t}.
\end{equation}
Likewise, if $\sigma\in S^{s_1\dots s_{t-1}l_t}$,
then, by the Filippov-Jacobi identity, $A_{\sigma}$ can be written as a sum of
one element $A_\rho$ with $\rho\in S^{l_1\dots l_t}$
as in (\ref{rho}), and
 $n-1$ elements $A_{\tau}$ with $\tau\in S^{l_1\dots l_{t-1}s_t}$, 
such that $\{\tau(1), \dots, \tau(n-2)\}\subset\{\sigma(n),\dots,
\sigma(2n-2)\}$. As above, there exists only one $i$
such that $\{\tau(1), \dots, \tau(n-2)\}\cup\{i\}=\{\sigma(n),\dots,
\sigma(2n-2)\}$. Then $i\leq 2n-t-1$, $i\neq \tau(1),\dots, \tau(n-2)$.
It follows that 
\begin{equation}
Q_{s_1\dots s_{t-1}l_t}=(n-t+1)
Q_{l_1\dots l_{t-1}s_t}+
(-1)^{n-1}Q_{l_1\dots l_t}.
\end{equation}
 Then,
using (\ref{l+s}) and the inductive hypotheses $Q_{s_1\dots s_{t-1}}=0=
Q_{l_1\dots l_{t-1}}$, we get the following
system of linear equations:
\begin{equation}
\left\{
\begin{array}{l}
Q_{s_1\dots s_t}=(n-t)Q_{l_1\dots l_t}\\
Q_{s_1\dots s_{t-1}l_t}=(n-t+1)Q_{l_1\dots l_{t-1}s_t}+(-1)^{n-1}Q_{l_1\dots l_t}\\
Q_{s_1\dots s_t}+Q_{s_1\dots s_{t-1}l_t}=0\\
Q_{l_1\dots l_t}+Q_{l_1\dots l_{t-1}s_t}=0,
\end{array}
\right.
\label{system1}
\end{equation}
whose determinant is equal to $(-1)^n+1$. It follows that if $n$ is even
then
$Q_{l_1\dots l_t}=0=Q_{s_1\dots s_t}=Q_{s_1\dots s_{t-1}l_t}=
Q_{l_1\dots l_{t-1}s_t}$, hence (\ref{thesis}) is proved.

Now assume that $n$ is odd. Then (\ref{system1}) reduces to
\begin{equation}\left\{
\begin{array}{l}
Q_{s_1\dots s_t}=(n-t)Q_{l_1\dots l_t}\\
Q_{s_1\dots s_{t-1}l_t}=(n-t+1)Q_{l_1\dots l_{t-1}s_t}+Q_{l_1\dots l_t}\\
Q_{s_1\dots s_t}+Q_{s_1\dots s_{t-1}l_t}=0.
\end{array}
\label{system1B}
\right.
\end{equation}

 Using the Filippov-Jacobi identity as above,
 one gets the following system
of linear equations: 
\begin{equation}
\left\{
\begin{array}{l}
Q_{s_1\dots s_{t-2}l_{t-1}l_t}=(n-t+2)Q_{l_1\dots l_{t-2}s_{t-1}s_t}
-Q_{l_1\dots l_{t-2}s_{t-1}l_t}
+Q_{l_1\dots l_{t-2}l_{t-1}s_{t}}\\
Q_{s_1\dots s_{t-2}l_{t-1}s_t}=(n-t+1)Q_{l_1\dots l_{t-2}s_{t-1}l_t}
+Q_{l_1\dots l_{t-2}l_{t-1}l_t}.
\end{array}
\right.
\label{system2}
\end{equation}
Besides, using the inductive hypotheses, we get:
\begin{equation}
\left\{
\begin{array}{l}
Q_{s_1\dots s_{t-2}l_{t-1}l_t}+Q_{s_1\dots s_{t-2}l_{t-1}s_t}=
Q_{s_1\dots s_{t-2}l_{t-1}}=0\\
Q_{l_1\dots l_{t-2}s_{t-1}s_t}+Q_{l_1\dots l_{t-2}s_{t-1}l_t}=
Q_{l_1\dots l_{t-2}s_{t-1}}=0\\
Q_{l_1\dots l_{t-2}l_{t-1}s_t}+Q_{l_1\dots l_{t-2}l_{t-1}l_t}=
Q_{l_1\dots l_{t-2}l_{t-1}}=0.
\end{array}
\right.
\label{system2B}
\end{equation}
Taking the sum of the two equations in (\ref{system2}),
 and using the three
equations in (\ref{system2B}), we get:
$Q_{l_1\dots l_{t-2}s_{t-1}s_t}=Q_{l_1\dots l_{t-2}s_{t-1}l_t}=0$.
By arguing in the same way, one shows that 
\begin{equation}
\label{eq:3.8}
Q_{l_1\dots l_ks_{k+1}l_{k+2}
\dots l_t}=0 ~~~\mbox{for every} ~~k=0,\dots t-2.
\end{equation}
Finally, using the Filippov-Jacobi identity as above, one
gets the following system of linear equations:
$$
\left\{
\begin{array}{l}
Q_{l_1\dots l_t}=nQ_{s_1\dots s_t}+(-1)^{n+t-2}Q_{s_1\dots s_{t-1}l_t}+
\dots -Q_{s_1l_2s_3\dots s_t}+Q_{l_1s_2\dots s_t}\\
Q_{s_1\dots s_{t-2}l_{t-1}s_t}=(n-t+1)Q_{l_1\dots l_{t-2}s_{t-1}l_t}+Q_{l_1\dots l_t}\\
\vdots\\
Q_{l_1s_2\dots s_{t}}=(n-t+1)Q_{s_1l_2\dots l_{t}}+Q_{l_1\dots l_t}
\end{array}\right.
$$
which reduces, by (\ref{eq:3.8}), to the following
\begin{equation}
\label{eq:3.9}
\left\{
\begin{array}{l}
Q_{l_1\dots l_t}=nQ_{s_1\dots s_t}+(-1)^{n+t-2}Q_{s_1\dots s_{t-1}l_t}+
\dots -Q_{s_1l_2s_3\dots s_t}+Q_{l_1s_2\dots s_t}\\
Q_{s_1\dots s_{t-2}l_{t-1}s_t}=Q_{l_1\dots l_t}\\
\vdots\\
Q_{l_1s_2\dots s_{t}}=Q_{l_1\dots l_t}.
\end{array}\right.
\end{equation}
System (\ref{eq:3.9}) implies the following equation:
$$Q_{l_1\dots l_t}=nQ_{s_1\dots s_t}+(-1)^{n+t-2}Q_{s_1\dots s_{t-1}l_t}+
\frac{(-1)^t+1}{2}Q_{l_1\dots l_t}.$$
This equation together with (\ref{system1B}) form a system of four linear 
equations in four indeterminates,
whose determinant is equal to $(n-t+1)((-1)^{n-t}(n-t)+\frac{1-(-1)^t}{2}
-n(n-t))$. It is different from 0 for every $t=1,\dots,n-2$. Hence
$Q_{l_1\dots l_t}=0=Q_{l_1\dots l_{t-1}s_t}=Q_{s_1\dots s_{t-1}l_t}=
Q_{s_1\dots s_t}$, and (\ref{thesis}) is proved.
\hfill$\Box$

\medskip

\noindent
{\bf Proof of Theorem \ref{th:0.2}.}  
Any element of the subalgebra generated by $L_{-1}$ and $\mu$ is a linear
combination of elements of the form:
\begin{equation}
[\dots[[\dots[[[\dots[\mu,a_1],\dots, a_s], \mu], b_1],\dots, b_k],\mu],\dots]
\label{general}
\end{equation}
with $a_1,\dots, a_s, b_1,\dots, b_k, \dots$ in $L_{-1}$.
By Lemma \ref{Askar}, every element of the form 
$[[[\dots[\mu,a_1],\dots, a_s], \mu]$ is either 0 or an 
element in $L_{-1}$, therefore we can assume
 that $\mu$ appears only once  in (\ref{general}), i.e., 
any element of $L$ lies in 
$[\dots[\mu, L_{-1}],\dots, L_{-1}]$.
\hfill$\Box$

\begin{remark}\em 
\label{rem:3.2}
We conjecture that Theorem \ref{th:0.2} holds also in non-zero characteristic 
if $\char \F \geq n$. Our argument works for $\char \F >(n-1)^2$.
The following example shows that 
Theorem \ref{th:0.2} (and Theorem \ref{th:0.1}) fails if $0< \char \F < n$.
Let $\fg = \F a$ be a 1-dimensional odd space, which we endow by the following
$n$-bracket: $[a,a,\ldots,a]=a$. The Filippov-Jacobi identity holds
if $n=sp+1$, where $p=\char \F$ and $s$ is a positive integer.   
However, $\Lie \g$ is not of the form described by Theorem \ref{th:0.2}.
\end{remark}

\section{Classification of irreducible admissible pairs}
\label{sec:4}
First, we briefly recall some examples of $\Z$-graded linearly compact Lie
superalgebras over a field $\F$ of characteristic 0, and some of their 
properties. For more details, see \cite{K2}
and \cite{CK3}.

Given a finite-dimensional vector superspace $V$ of dimension $(m|n)$
(i.e.\ $\dim V_{\0}=m$, $\dim V_{\1}=n$), the universal Lie superalgebra
$W(V)$ is isomorphic to the Lie superalgebra $W(m,n)$ of continuous derivations
of the tensor product $\F(m,n)$ of the algebra of formal power series in $m$
commuting variables $x_1,\dots, x_m$ and the Grassmann algebra in $n$
anti-commuting variables $\xi_1,\dots,\xi_n$. Elements of $W(m,n)$ can be 
viewed as linear differential operators of the form
$$X=\sum_{i=1}^m P_i(x,\xi)\frac{\partial}{\partial x_i}+
\sum_{j=1}^n Q_j(x,\xi)\frac{\partial}{\partial \xi_j}, ~P_i, Q_j\in\F(m,n).$$
The Lie superalgebra $W(m,n)$ is simple linearly compact (and
it is finite-dimensional if and only if $m=0$).

Letting $\deg x_i=-\deg\frac{\partial}{\partial x_i}=k_i$,
$\deg\xi_i=-\deg\frac{\partial}{\partial \xi_i}=s_i$, where
$k_i, s_i\in\Z$, defines a $\Z$-grading on $W(m,n)$, called the $\Z$-grading
of type $(k_1,\dots,k_m|s_1,\dots,s_n)$. Any $\Z$-grading of $W(m,n)$ is
conjugate (i.e.\ can be mapped by an automorphism of $W(m,n)$)
to one of these. Clearly, such a grading has finite
depth $d$ (meaning that  $W(m,n)_j\neq 0$ if and only if $j\geq -d$) if and
only if $k_i\geq 0$ for all $i$. It is easy to show that the depth $d=1$ if
all $k_i$'s and $s_i$'s are 0 or 1, or if all $k_i$'s are 0, 
$s_j=-1$ for some $j$, and $s_i=0$ for every $i\neq j$. 

Now we shall describe some closed (hence linearly compact) subalgebras
of $W(m,n)$.

First, given a subalgebra $L$ of $W(m,n)$, a continuous linear map $Div: L
\rightarrow \F(m,n)$ is called a divergence if the action
$\pi_{\lambda}$ of $L$ on $\F(m,n)$, given by
$$\pi_{\lambda}(X)f=Xf+(-1)^{p(X)p(f)}\lambda f Div X, \,\, X \in L,$$
is a representation of $L$ in $\F(m,n)$ for any $\lambda\in\F$. Note that
$$S'_{Div}(L):=\{X\in L~|~ Div X=0\}$$
is a closed subalgebra of $L$. We denote by $S_{Div}(L)$ its derived subalgebra
(recall that the derived subalgebra of $\g$ is $[\g,\g]$). An example
of a divergence on $L=W(m,n)$ is the following, denoted by $div$:
$$div(\sum_{i=1}^m P_i\frac{\partial}{\partial x_i}+
\sum_{j=1}^n Q_j\frac{\partial}{\partial \xi_j})=
\sum_{i=1}^m \frac{\partial P_i}{\partial x_i}+
\sum_{j=1}^n (-1)^{p(Q_j)} \frac{\partial Q_j}{\partial \xi_j}.$$
Hence for any $\lambda\in\F$ we get the representation $\pi_\lambda$ of
$W(m,n)$ in $\F(m,n)$. Also, we get closed subalgebras $S'_{div}(W(m,n))
\supset S_{div}(W(m,n))$ denoted by $S'(m,n)\supset S(m,n)$. Recall
that $S'(m,n)=S(m,n)$ is simple if $m>1$, and
\begin{equation}
\label{eq:4.1}
S'(1,n)=S(1,n)\oplus\F\xi_1\dots\xi_n\frac{\partial}{\partial x_1},
\end{equation}
where $S(1,n)$ is a simple ideal.

The $\Z$-gradings of type $(k_1,\dots,k_m|s_1,\dots,s_n)$ of $W(m,n)$
induce ones on $S'(m,n)$ and $S(m,n)$ and any $\Z$-grading is conjugate to 
those. The description of $\Z$-gradings of depth 1 for $S'(m,n)$ and $S(m,n)$
is the same as for $W(m,n)$.

Next examples of subalgebras of $W(m,n)$, needed in this paper, are of the
form
$$L(\omega)=\{X\in W(m,n)~|~X\omega=0\},$$
where $\omega$ is a differential form.

In the case $m=2k$ is even, consider the symplectic differential form
$$\omega_s=2\sum_{i=1}^kdx_i\wedge dx_{k+i}+\sum_{i=1}^n
d\xi_id\xi_{k-i+1}.$$
The corresponding subalgebra $L(\omega_s)$ is denoted by $H'(m,n)$ and is
called a Hamiltonian superalgebra. This superalgebra is simple,
hence coincides with its derived subalgebra $H(m,n)$, unless $m=0$, when the
Hamiltonian superalgebra is finite-dimensional.

It is convenient to consider the ``Poisson'' realization of $H(m,n)$. For that
let $p_i=x_i$, $q_i=x_{k+i}$, $i=1,\dots,k$, and introduce on $\F(m,n)$
the structure of a Poisson superalgebra $P(m,n)$ by letting the non-zero
brackets between generators to be as follows:
$$\{p_i,q_i\}=1=\{\xi_i,\xi_{n-i+1}\},$$
and extend by the Leibniz rule. Then the map $P(m,n)\rightarrow H'(m,n)$,
given by $f\mapsto \sum_{i=1}^k(\frac{\partial f}{\partial p_i}
\frac{\partial}{\partial q_i}- \frac{\partial f}{\partial q_i}
\frac{\partial}{\partial p_i})-(-1)^{p(f)}
\sum_{i=1}^k\frac{\partial f}{\partial \xi_i}
\frac{\partial}{\partial \xi_{k-i+1}},$
defines a surjective Lie superalgebra homomorphism with kernel $\F 1$. Thus,
$H'(m,n)=P(m,n)/\F 1$. In this realization 
$H(0,n)$ is spanned by all monomials in $\xi_i$ mod $\F 1$ except for the one of top degree, and we have:
\begin{equation}
\label{eq:4.2}
H'(0,n)=H(0,n)\oplus\F\xi_1\dots\xi_n.
\end{equation}
Note that $H(0,n)$ is simple if and only if $n\geq 4$.

All $\Z$-gradings of depth 1 of $H'(0,n)$ are, up to conjugacy, those of type
$(|1,\dots,1)$, $(|1,0,\dots,0,$ $-1)$,
and $(|\underbrace{1,\dots,1}_{n/2},0,\dots,0)$, if $n$ is even \cite{CK4}.

Another example is $HO(n,n)=L(\omega_{os})\subset W(n,n)$, where
$\omega_{os}=\sum_{i=1}^n dx_id\xi_i$ is an odd symplectic form. This Lie
superalgebra is simple if and only if $n\geq 2$. It contains
the important for this paper subalgebra $SHO'(n,n)=HO(n,n)\cap S'(n,n)$.
Its derived subalgebra $SHO(n,n)$ is simple if and only if $n\geq 3$.

Again, it is convenient to consider a ``Poisson'' realization of 
$HO(n,n)$. For this
consider the Buttin bracket on $\Pi\F(n,n)$:
$$\{f,g\}_B=\sum_{i=1}^n(\frac{\partial f}{\partial x_i}
\frac{\partial g}{\partial \xi_i}- (-1)^{p(f)}
\frac{\partial f}{\partial \xi_i}\frac{\partial g}{\partial x_i}).$$
This is a Lie superalgebra, which we denote by $PO(n,n)$, and the map
$PO(n,n)\rightarrow HO(n,n)$, given by
$$f\mapsto \sum_{i=1}^n(\frac{\partial f}{\partial x_i}
\frac{\partial }{\partial \xi_i}- (-1)^{p(f)}
\frac{\partial f}{\partial \xi_i}\frac{\partial}{\partial x_i})$$
is a surjective Lie superalgebra homomorphism, whose kernel is $\F 1$. Thus, 
$HO(n,n)=P(n,n)/\F 1$. In this realization we have:
$$SHO'(n,n)=\{f\in P(n,n)/\F 1 ~|~ \Delta f=0\},$$
where $\Delta=\sum_{i=1}^n\frac{\partial}{\partial x_i}
\frac{\partial}{\partial \xi_i}$ is the odd Laplace operator. Then $SHO(n,n)$ 
is an ideal of codimension 1 in $SHO'(n,n)$, and we have:
\begin{equation}
\label{eq:4.3}
SHO'(n,n)=SHO(n,n)\oplus\F\xi_1\dots\xi_n.
\end{equation}
All $\Z$-gradings of depth 1 of $SHO'(n,n)$ are, up to conjugacy, those of type
$(1,\dots,1|1,\dots,1)$, $(0,\dots,0,1|0,\dots,0,-1)$, and
$(\underbrace{1,\dots,1}_k,0,\dots,0|\underbrace{0,\dots,0}_k,1,\dots,1)$,
where $k=0,\dots,n$ \cite{CK4}.

The next important for us example is  
$$KO(n,n+1)=\{X\in W(n,n+1)~|~X\omega_{oc}=f\omega_{oc}~~\mbox{for some}~
f\in\F(n,n+1)\},$$
where $\omega_{oc}=d\xi_{n+1}+\sum_{i=1}^n(\xi_idx_i+x_id\xi_i)$ is an odd
contact form. This superalgebra is simple for all $n\geq 1$. Another 
realization of this Lie superalgebra is $PO(n,n+1)=\Pi\F(n,n+1)$ with the
bracket
 $\{f,g\}_{BO}=(2-E)f\frac{\partial g}{\partial\xi_{n+1}}-
(-1)^{p(f)}\frac{\partial f}{\partial\xi_{n+1}}(2-E)g-
\sum_{i=1}^n(\frac{\partial f}{\partial x_i}\frac{\partial g}{\partial\xi_{i}}
-(-1)^{p(f)}\frac{\partial f}{\partial\xi_{i}}\frac{\partial g}{\partial x_i})$,
 where $E=\sum_{i=1}^n(x_i\frac{\partial}{\partial x_i}+
\xi_i\frac{\partial}{\partial\xi_i})$. The isomorphism $PO(n,n+1)\rightarrow
KO(n,n+1)$ is given by $f\mapsto (2-E)f\frac{\partial}{\partial\xi_{n+1}}+
(-1)^{p(f)}\frac{\partial f}{\partial\xi_{n+1}}E-
\sum_{i=1}^n(\frac{\partial f}{\partial x_i}\frac{\partial}{\partial\xi_{i}}
-(-1)^{p(f)}\frac{\partial f}{\partial\xi_{i}}\frac{\partial}{\partial x_i}).
$
It turns out that for each $\beta\in\F$ the Lie superalgebra $KO(n,n+1)$
admits a divergence
$$div_{\beta}f=\Delta f+(E-n\beta)\frac{\partial f}{\partial\xi_{n+1}}, \,\,
f\in PO(n,n+1).$$
We let
$$SKO'(n,n+1;\beta)=\{f\in PO(n,n+1)~|~ div_\beta f=0\}.$$
This Lie superalgebra is not always simple, but its derived algebra, denoted by
$SKO(n,n+1;\beta)$, is simple if and only if $n\geq 2$. In fact,
$SKO'(n,n+1;\beta)=SKO(n,n+1;\beta)$, unless $\beta=1$ or $\beta=\frac{n-2}{n}$.
In the latter cases $SKO(n,n+1;\beta)$ is an ideal of codimension 1
in $SKO'(n,n+1;\beta)$, and we have:
\begin{equation}
\label{eq:4.4}
SKO'(n,n+1;1)=SKO(n,n+1;1)+\F\xi_1\dots\xi_{n+1},
\end{equation}
\begin{equation}
\label{eq:4.5}
SKO'(n,n+1;\frac{n-2}{n})=SKO(n,n+1;\frac{n-2}{n})+\F\xi_1\dots\xi_{n}.
\end{equation}
All $\Z$-gradings of depth 1 of $SKO'(n,n+1;\beta)$ are, up to conjugacy, of 
type $(0,\dots,0,1|0,\dots,0,-1,0)$, and
$(\underbrace{1,\dots,1}_k,0,\dots,0|\underbrace{0,\dots,0}_k,1,
\dots,1)$, where $k=0,\dots,n$ \cite{CK4}.

\begin{theorem}
\label{thm:4.1}
Let $(L=\oplus_{j=-1}^{n-1}L_j,\mu)$ be an irreducible admissible pair over
an algebraically closed field $\F$ of characteristic 0, where $L$ is a
linearly compact Lie superalgebra, and $n>2$. Then
\begin{itemize}
\item[(a)] $L=\oplus_{j=-1}^{n-1}L_j$ is a semidirect product of the simple
ideal $S=\oplus_{j=-1}^{n-2}L_j$ and the 1-dimensional subalgebra 
$L_{n-1}=\F\mu$, where $\mu$ is an outer derivation of $L$, such that
$[\mu,L_0]=0$.
\item[(b)] The pair $(L,\mu)$ is isomorphic to one of the following four
irreducible admissible pairs:
\begin{itemize}
\item[(i)] $(H'(0,n+1), \xi_1\dots\xi_{n+1})$, $n\geq 3$, with the
grading of type $(|1,\dots,1)$;
\item[(ii)] $(SHO'(n,n), \xi_1\dots\xi_{n})$, $n\geq 3$, with the
grading of type $(0,\dots,0|1,\dots,1)$;
\item[(iii)] $(SKO'(n-1,n;1), \xi_1\dots\xi_{n-1}\xi_{n})$, $n\geq 3$, with the
grading of type $(0,\dots,0|1,\dots,1)$;
\item[(iv)] $(S(1,n-1), \xi_1\dots\xi_{n-1}\frac{\partial}{\partial x})$, $n\geq 3$, with the
grading of type $(0|1,\dots,1)$.
\end{itemize}
\end{itemize}
\end{theorem}
{\bf Proof.} The decomposition $L=S\rtimes\F\mu$ in $(a)$ follows from Theorem 
\ref{th:0.2}. The fact that $S$ is simple is proved in the same way as in
\cite[Theorem 7.3]{CK4}. Indeed, $S$ is the minimal among non-zero closed ideals of $L$, since if $I$ is a  non-zero closed ideal of $L$, 
then $I\cap L_{-1}\neq\emptyset$ by transitivity, hence,
by irreducibility, $I\cap L_{-1}=L_{-1}$, from which it follows that $I$ contains $S$.
Next, by the super-analogue of 
Cartan-Guillemin's theorem
\cite{B,G1}, established in \cite{FK}, 
$S=S'\hat{\otimes}\Lambda(m,h)$, for some simple linearly compact
Lie superalgebra $S'$ and some
$m,h\in\Z_{\geq 0}$, and $\mu$ lies in $Der(S'\hat{\otimes}\O(m,h))$. Since
$Der(S'\hat{\otimes}\O(m,h))=Der S'\hat{\otimes}\O(m,h)+1\otimes W(m,h)$
\cite{FK}, we have: 
$\mu=\sum_i(d_i\otimes a_i) +1\otimes \mu'$ for some $d_i\in Der S'$,
$a_i\in\O(m,h)$ and $\mu'\in W(m,h)$.

First consider the case when $\mu$ is even.
Then $\mu'$ is an even element of $W(m,h)$ hence, by the minimality of
the ideal $S'\hat{\otimes} \O(m,h)$, $h=0$. Now suppose 
$m\geq 1$. If $\mu'$ lies
in the non-negative part of $W(m,0)$ with the grading of type $(1,\dots,1|)$, then the
ideal generated by $S'x_1$ is a proper $\mu$-invariant ideal of 
$S'\hat{\otimes} \O(m,0)$,
contradicting its minimality. Therefore we may assume, up to a linear change 
of indeterminates, that $\mu'=\frac{\partial}{\partial x_1}+D$, for some
derivation $D$ lying in the non-negative part of
$W(m,0)$. Since $\mu$ lies in $L_{n-1}$, we have $\deg(x_1)=-n+1$, but this
is a contradiction since the $\Z$-grading of $L$ has depth 1.
It follows that $m=0$.  

Now consider the case when $\mu$ is odd.
Consider the grading of $W(m,h)$ of type $(1,\dots,1|1,\dots,1)$, and
denote by $W(m,h)_{\geq 0}$ its non-negative part. If $\mu'\in
W(m,h)_{\geq 0}$, then the minimality of the ideal
$S'\hat{\otimes}{\cal O}(m,h)$ implies $m=h=0$.
Now suppose that $h\geq 1$
and that $\mu'$ has a non-zero projection on $W(m,h)_{-1}$. Then,  up
to a linear change of indeterminates, 
$\mu'=\frac{\partial}{\partial\xi_1}+D$ for some odd derivation $D\in
W(m,h)_{\geq 0}$.   
Since $\mu$ lies in $L_{n-1}$, we have $\deg(\xi_1)=-n+1 < -1$.
Since $L\supset S'\xi_1$ and the grading of $L$ has depth 1, it follows that 
every element in $S'$ has  positive degree, but this is a contradiction since 
$S'$ is simple. This concludes the proof of the simplicity of $S$.


In order to prove $(b)$, note that the grading operator $D$ of the simple
$\Z$-graded Lie superalgebra $S$ is an outer derivation (since $[\mu,L_0]=0, 
D\notin L_0$). Another outer derivation of $S$ is $\mu$. From the
classification of simple linearly compact Lie superalgebras
\cite{K1}, \cite{K2} and their derivations in \cite{K1}, \cite{K2},
\cite[Proposition 1.8]{CK1} (see also Lemma \ref{lemma:4.3} below), 
we see that the only possibilities for $L$
are $H'(0,n+1)$, $SHO'(n,n)$, $SKO'(n-1,n;1)$,
$SKO'(n-1,n;\frac{n-2}{n})$, and $S'(1,n)$ for $n\geq 3$. From the
description of $\Z$-gradings of depth 1 of these Lie superalgebras, 
given above,
it follows that $L=SKO'(n-1,n;\frac{n-1}{n})$ is ruled out, whereas for the
remaining four possibilities for $L$ only the grading of type
$(0,\dots,0|1,\dots,1)$ is possible, and for them indeed $L_{n-1}=\F\mu$,
where $\mu$ is as described above. It is immediate to check that in these four
cases the pair $(L,\mu)$ is admissible. Irreducibility of the 
$L_0$-module $L_{-1}$ follows
automatically from the simplicity of $S$ since its depth is 1.
 \hfill$\Box$

\begin{remark}\em
\label{Rk:4.2}
In cases $(i)$--$(iv)$ of Theorem \ref{thm:4.1}$(b)$ the subalgebra
$L_0$ and the $L_0$-module $\Pi L_{-1}$ are as follows:
\begin{itemize}
\item[$(i)$] $L_0\cong so_{n+1}(\F)$, $\Pi L_{-1}=\F^{n+1}$ with the
standard action of $so_{n+1}(\F)$;
\item[$(ii)$] $L_0\cong S(n,0)$, $\Pi L_{-1}=\F[[x_1,\dots,x_n]]/\F 1$,
 where $\F[[x_1,\dots,x_n]]$ is the
standard module over $S(n,0)$;
\item[$(iii)$] $L_0\cong W(n-1,0)$, $\Pi L_{-1}=\F[[x_1,\dots,x_{n-1}]]$,
which carries the representation $\pi_{\lambda=-1}$ of
$W(n-1,0)$;
\item[$(iv)$] $L_0\cong W(1,0)\ltimes sl_{n-1}(\F[[x]])$, 
$\Pi L_{-1}=\F^{n-1}\otimes\F[[x]]$ with the standard action of 
$sl_{n-1}(\F[[x]])$ and the representation $\pi_{\lambda=-1/(n-1)}$
of  $W(1,0)$ on $\F[[x]]$.
\end{itemize}
\end{remark}

As we have seen, an important part of the classification of irreducible
admissible pairs is the description of derivations of simple linearly
compact Lie algebras. This description is based on the following simple
lemma.
\begin{lemma}
\label{lemma:4.3}
Let $L$ be a linearly compact Lie superalgebra and let $\mathfrak{a}$ be a 
reductive subalgebra of $L$ (i.e. the adjoint representation of 
$\mathfrak{a}$ on $L$ decomposes in a direct product of finite-dimensional
irreducible $\mathfrak{a}$-modules). Then any continuous derivation of $L$
is a sum of an inner derivation and a derivation commuting with the adjoint
action of $\mathfrak{a}$.
\end{lemma}
{\bf Proof.} \cite{K1} We have closed $\mathfrak{a}$-submodules:
$$\Inder L\subset \Der L\subset End L\,,$$
where $\Inder L$ and $\Der L$ denote the subspaces of all inner derivations and
all continuous derivations of the Lie superalgebra $L$
in the space of continuous endomorphisms of the linearly compact vector space
$L$. Since $L=\prod_j V_j$, where $V_j$ are finite-dimensional irreducible
$\mathfrak{a}$-modules, we have:
$\End L=\prod_{i,j} \Hom(V_i, V_j)$, hence $\End L$, and therefore $\Der L$, 
decomposes into a direct product of irreducible $\mathfrak{a}$-submodules.
Hence
$$\Der L=\Inder L\oplus V,$$
where $V$ is an $\mathfrak{a}$-submodule. But $\mathfrak{a}V\subset \Inder L$
since $\Inder L$ is an ideal in $\Der L$. Hence $\mathfrak{a}V=0$, i.e., 
any derivation from $V$ commutes with the adjoint action of $\mathfrak{a}$
on $L$.
\hfill$\Box$

\section{Classification of simple linearly compact $n$-Lie algebras over a field
of characteristic 0, and their derivations}
\label{sec:5}

{\bf Proof of Theorem \ref{th:0.1}} By Proposition \ref{prop:2.4},
the classification of simple linearly compact $n$-Lie algebras is equivalent to
the classification of admissible pairs $(L,\mu)$, for which $L$ is linearly
compact. The list of the latter consists of the four examples $(i)$--$(iv)$
given in Theorem \ref{thm:4.1}$(b)$. It is easy to see that the corresponding
$n$-Lie algebras are $O^n$, $S^n$, $W^n$ and $SW^n$. (By Lemma 
\ref{lem:1.4}$(a)$, we automatically get from $[\mu,L_0]=0$ that 
the Filippov-Jacobi identity indeed holds.) \hfill$\Box$

\medskip

The notation for the four simple $n$-Lie algebras comes from the following
fact.
\begin{proposition}
\label{prop:5.1}
\begin{itemize}
\item[(a)] The Lie algebra of continuous derivations of the $n$-Lie algebras
$O^n$, $S^n$, $W^n$ and $SW^n$ is isomorphic to $so_{n+1}(\F)$, $S(n,0)$,
$W(n-1,0)$ and $W(1,0)\ltimes sl_{n-1}(\F[[x]])$, respectively. Its 
representation on the $n$-Lie algebra is described in Remark \ref{Rk:4.2}.
\item[(b)] All continuous derivations of a simple linearly compact $n$-Lie
algebra $\fg$ over an algebraically closed field of characteristic 0 
lie in the closure $\Inder \fg$ of the span of the inner ones.
\end{itemize}
\end{proposition}
{\bf Proof.} Let $\g$ be one of the four simple $n$-Lie algebras and let
$\Der \g$ be the Lie algebra of all continuous  derivations of $\g$. Then
$L_0 :=\Inder \fg$ is an ideal of $\Der \g$.
By Remark \ref{Rk:4.2}, $L_0$ is isomorphic to the Lie algebras listed in 
$(a)$. But all derivations of the Lie algebras 
$L_0=so_{n+1}(\F)$, $W(n-1,0)$ and 
$W(1,0)\ltimes sl_{n-1}(\F[[x]])$ are inner, and 
$\Der S(n,0)=S(n,0)\oplus \F E$,
where $E=\sum_i x_i\frac{\partial}{\partial x_i}$. This is well known,
except for the case $L=W(1,0)\ltimes sl_{n-1}(\F[[x]])$.
We apply Lemma \ref{lemma:4.3} to this case, taking 
$\mathfrak{a}=\F x\frac{d}{dx}\oplus sl_{n-1}(\F)$. If $D$ is an endomorphism 
of the vector space $L$, commuting with $\mathfrak{a}$, we have, by Schur's lemma
:
$$D(x^k\frac{d}{dx})=\alpha_kx^k\frac{d}{dx}, ~~D(x^ka)=\beta_kx^ka,
~{\mbox{for}}~ a\in sl_{n-1}(\F), {\mbox{where}}\,\, 
\alpha_k , \beta_k \in \F.$$
Since $D$ is also a derivation of $L$, we conclude that $D$ is a 
multiple of $\ad x\frac{d}{dx}$.

Let now $D\in \Der \g\setminus(\Inder \g=L_0)$. Since $[D,L_0]\subset L_0$,
$D$ induces a derivation of $L_0$. Since all derivations of $L_0$ are inner,
except for $E$ in the case $\g\cong S^n$, but $E$ is not a derivation of
$\g$, we conclude that there exists $a\in L_0$, such that 
$D|_{L_0}=\ad a|_{L_0}$. Therefore $D'=D-a$ commutes with the action of $L_0$
on $L_{-1}$. But the latter representation is described in Remark \ref{Rk:4.2},
and, clearly, in all cases the only operators, commuting with the 
representation operators of $L_0$
on $L_{-1}$, are scalars. Since a non-zero scalar cannot be a derivation of
$\g$, we conclude that $D'=0$, hence $\Der \g=L_0$. \hfill$\Box$

\medskip

In conclusion we discuss $\F$-forms of the four simple $n$-Lie algebras, where
$\F$ is a field of characteristic 0. Let $\overline{F}\supset \F$ be the
algebraic closure of $\F$. Given a linearly compact $n$-Lie algebra
$\g$ over $\overline{\F}$, its $\F$-form is defined as 
an $n$-Lie algebra $\g^{\F}$
over $\F$, such that $\overline{F}\otimes_{\F}\g^{\F}$ is isomorphic to $\g$.

Due to the bijection given by Proposition \ref{prop:2.4}, the $\F$-forms
of $\g$ are in one-to-one correspondence with the $\F$-forms of the 
$\Z$-graded Lie superalgebras $S_{\g}=[\Lie \g, \Lie \g]$. But the latter
are parameterized by the set $H^1(Gal, Aut S_{\g})$, where $Gal$ is the
Galois group of $\overline{\F}$ over $\F$, and $Aut S_{\g}$ is the group 
of continuous automorphisms of the Lie superalgebra $S_{\g}$, preserving
its $\Z$-grading (cf.\ \cite{CK2}).

By the method of \cite{CK2} it is easy to compute the group $Aut S_{\g}$,
using Remark \ref{Rk:4.2}.

\begin{proposition}
\label{prop:5.2}
One has:
$$Aut S_{\g}=G_{\g}\ltimes{\cal U},$$
where ${\cal U}$ is a prounipotent group and $G_{\g}$ is a reductive
group, isomorphic to $O_{n+1}(\overline{\F})$,
$GL_{n}(\overline{\F})$, $GL_{n-1}(\overline{\F})$ and $\overline{\F}^{\times}
\times SL_{n-1}(\overline{F})$, if $\g$ is isomorphic to $O^n$, $S^n$,
$W^n$ and $SW^n$ over $\overline{\F}$, respectively.
\end{proposition}
We have $H^1(Gal, Aut S_{\g})=H^1(G_{\g}, Gal)$ (see, e.g., \cite{CK2}).
Furthermore, $H^1(G_{\g}, Gal)=1$ in the last three cases of
Proposition \ref{prop:5.2}, hence the only $\F$-forms of $S^n$, $W^n$ and 
$SW^n$ over $\overline{\F}$ are $S^n$, $W^n$ and $SW^n$ over $\F$.
Finally, it follows from \cite{CK2} that the $\F$-forms of the $\Z$-graded
Lie superalgebra $H(0,n+1)$ are the derived algebras of the Lie superalgebras
$P/\F 1$, where $P$ is a Poisson algebra, defined by (\ref{eq:0.6}).
Hence $\F$-forms of $O^n$ are vector product $n$-Lie algebras on 
$\F^{n+1}$, $n\geq 3$, with a non-degenerate symmetric bilinear
form (up to isomorphism, these $n$-Lie algebras depend on the
equivalence class of the bilinear form up to a non-zero factor).

\section*{Appendix A}
Below we list all known examples of infinite-dimensional simple $n$-Lie
algebras over an algebraically closed field $\FF$ of characteristic 0
for $n \geq 3$.

Let $A$ be a commutative associative algebra over $\FF$ and let
$\fg$ be a Lie algebra of derivations of $A$, such that $A$
contains no non-trivial $\fg$-invariant ideals.

\noindent{\bf Example 1.} \quad $S (A,\fg) =A$, where $\fg$ is an
$n$-dimensional Lie algebra with basis $D_1 ,\ldots ,D_n$, the
$n$-ary Lie bracket being
\begin{displaymath}
  [f_1,\ldots ,f_n] = \det \left(
    \begin{array}{lcl}
D_1 (f_1) & \ldots & D_1 (f_n)\\
    \hdotsfor{3}\\
D_n (f_1) & \ldots & D_1 (f_n)
    \end{array} \right) \, .
\end{displaymath}

\noindent{\bf Example 2.}  \quad $W (A,\fg)$, where $\fg$ is an
$n-1$-dimensional Lie algebra with basis $D_1 ,\ldots , D_{n-1}$,
the $n$-ary Lie bracket being 
\begin{displaymath}
  [f_1 ,\ldots ,f_n] = \det \left(
    \begin{array}{ccc}
      f_1 & \ldots & f_n\\
      D_1 (f_1) & \ldots & D_1 (f_n)\\
          \hdotsfor{3}\\
          D_{n-1} (f_1) & \ldots & D_{n-1}(f_n )
    \end{array}\right) \, .
\end{displaymath}

\noindent{\bf Example 3.} \quad $SW (A,D) =A^{\langle 1 \rangle}\oplus \cdots
\oplus A^{\langle n-1 \rangle}$ is the sum of $n-1$ copies of $A$ and $\fg =
\FF D$, the $n$-ary Lie bracket being the following.
For $h\in A$,  denote by $h^{\langle k \rangle}$ the corresponding element 
in $A^{\langle k \rangle}$, then
$$[f_1^{\langle j_1 \rangle},\dots, f_{n}^{\langle j_n \rangle}]=0, 
~~\mbox{unless}~ \{j_1, \dots, j_{n}\}
\supset \{1, \dots, n-1\};$$
$$[f_1^{\langle 1 \rangle},\dots, f_{k-1}^{\langle k-1 \rangle}, 
f_k^{\langle k \rangle}, f_{k+1}^{\langle k \rangle},
f_{k+2}^{\langle k +1 \rangle},\dots, f_{n}^{\langle n-1 \rangle}]=$$
$$(-1)^{k+n-1}(f_1\dots f_{k-1}
(D(f_k)f_{k+1}-f_kD(f_{k+1}))f_{k+2}\dots f_{n})^{\langle k \rangle},$$
extended on $SW(A,D)$ by anticommutativity.

\bigskip

It is an open problem whether there exist any other simple infinite-dimensional
$n$-Lie (super)algebras over an algebraically closed field of characteristic 0
if $n>2$. In particular are there any examples of infinite-dimensional simple $n$-Lie superalgebras
over a field of characteristic 0, which are not $n$-Lie algebras, if 
$n > 2$?  

\end{document}